\documentclass{article}
\usepackage{amsmath,amssymb,amsthm,amscd,dsfont}
\numberwithin{equation}{section} \allowdisplaybreaks
\begin{document}
\newtheorem{theorem}{Theorem}[section]
\newtheorem{defin}{Definition}[section]
\newtheorem{prop}{Proposition}[section]
\newtheorem{corol}{Corollary}[section]
\newtheorem{lemma}{Lemma}[section]
\newtheorem{rem}{Remark}[section]
\newtheorem{example}{Example}[section]
\title{A construction of Courant algebroids on foliated manifolds}
\author{{\small by}\vspace{2mm}\\Izu Vaisman}
\date{}
\maketitle
{\def\thefootnote{*}\footnotetext[1]%
{{\it 2000 Mathematics Subject Classification: 53C12, 53D17} .
}}
\begin{center} \begin{minipage}{12cm}
A{\footnotesize BSTRACT. For any transversal-Courant algebroid $E$
on a foliated manifold $(M,\mathcal{F})$, and for any choice of a
decomposition $TM=T\mathcal{F}\oplus Q$, we construct a Courant
algebroid structure on $T\mathcal{F}\oplus T^*\mathcal{F}\oplus E$.}
\end{minipage}
\end{center}
\section{Preliminaries}
General Courant algebroids were studied first in a paper by Liu,
Weinstein and Xu \cite{LWX}, which appeared in 1997 and became the
object of an intensive research since then. Courant algebroids
provide the framework for Dirac structures and generalized
Hamiltonian formalisms. In \cite{V1} we have introduced the notions
of transversal-Courant and foliated Courant algebroid, thereby
extending the framework to bases that are a space of leaves of a
foliation rather than a manifold. In the present note we show that a
transversal-Courant algebroid over a foliated manifold can be
extended to a foliated Courant algebroid. A similar construction for
Lie algebroids (which is a simpler case) was given in \cite{V1}. We
assume that the reader has access to the paper \cite{V1}, from which
we also take the notation, and he will consult \cite{V1} for the
various definitions and results that we use here. In this paper we
assume that all the manifolds, foliations, mappings, bundles, etc.,
are $C^\infty$-differentiable.

A Courant algebroid over the manifold $M$ is a vector bundle
$E\rightarrow M$ endowed with a symmetric, non degenerate, inner
product $g_E\in\Gamma\odot^2E^*$, with a bundle morphism
$\sharp_E:E\rightarrow TM$ called the {\it anchor} and a
skew-symmetric bracket $[\,,\,]_E:\Gamma E\times\Gamma
E\rightarrow\Gamma E$, such that the following conditions (axioms)
are satisfied:

1) $\sharp_E[e_1,e_2]_E=[\sharp_Ee_1,\sharp_Ee_2]$,

2) $im(\sharp_{g_E}\circ^t\sharp_E)\subseteq ker\,\sharp_E$,

3)
$\sum_{Cycl}[[e_1,e_2]_E,e_3]_E=(1/3)\partial_E\sum_{Cycl}g_E([e_1,e_2]_E,e_3)$,
$\partial_E=(1/2)\sharp_{g_E}\circ\,^t\sharp_E:T^*M\rightarrow E$,
$\partial_E f=\partial_E(df)$,

4)
$[e_1,fe_2]_E=f[e_1,e_2]_E+(\sharp_Ee_1(f))e_2-g(e_1,e_2)\partial_E
f $,

5) $(\sharp_Ee)(g_E(e_1,e_2))=g_E([e,e_1]_E+\partial_E g(e,e_1)
,e_2)+g_E(e_1,[e,e_2]_E+\partial_E g(e,e_2)).$\\ In these
conditions, $e,e_1,e_2,e_3\in\Gamma E$, $f\in C^\infty(M)$ and $t$
denotes transposition. Notice also that the definition of
$\partial_E$ is equivalent with the formula
\begin{equation}\label{partialcug}
g_E(e,\partial_Ef)=\frac{1}{2}\sharp_Ee(f).\end{equation}

The index $E$ will be omitted if no confusion is possible.

The basic example of a Courant algebroid was studied in \cite{C}
and it consists of the {\it big tangent bundle} $T^{big}M=TM\oplus
T^*M$, with the anchor $\sharp(X\oplus\alpha)=X$ and with
\begin{equation}\label{gC}
g(X_1\oplus\alpha_1,X_2\oplus\alpha_2)=\frac{1}{2}(\alpha_1(X_2)+\alpha_2(X_1)),
\end{equation}
\begin{equation}\label{crosetCou}
[X_1\oplus\alpha_1,X_2\oplus\alpha_2]=[X_1,X_2]\oplus
(L_{X_1}\alpha_2- L_{X_2}\alpha_1+\frac{1}{2}d(\alpha_1(X_2)
-\alpha_2(X_1))).
\end{equation}
(The notation $X\oplus \alpha$ instead of the
accurate $X+\alpha$ or $(X,\alpha)$ has the advantage of showing
the place of the terms while avoiding some of the parentheses. The
unindexed bracket of vector fields is the usual Lie bracket.)

Furthermore, let $\mathcal{F}$ be a foliation of the manifold $M$.
We denote the tangent bundle $T\mathcal{F}$ by $F$ and define the
transversal bundle $\nu\mathcal{F}$ by the exact sequence $$
0\rightarrow F
\stackrel{\iota}{\rightarrow}TM\stackrel{\psi}{\rightarrow}\nu\mathcal{F}
\rightarrow0,$$ where $\iota$ is the inclusion and $\psi$ is the
natural projection. We also fix a decomposition
\begin{equation}\label{descTM} TM=F\oplus
Q,\;Q=im(\varphi:\nu\mathcal{F}
\rightarrow TM),\,\psi\circ\varphi=id.,\end{equation} which implies
\begin{equation}\label{descT*M} T^*M=Q^*\oplus
F^*,\;Q^*=ann\,F,\,F^*=ann\,Q\approx T^*M/ann\,F,\end{equation}
where the last isomorphism is induced by the transposed mapping
$^t\iota$. The decompositions (\ref{descTM}), (\ref{descT*M})
produce a bigrading $(p,q)$ of the Grassmann algebra bundles of
multivector fields and exterior forms where $p$ is the $Q$-degree
and $q$ is the $F$-degree \cite{V0}.

The vector bundle $T^{big}\mathcal{F}=F\oplus (T^*M/ann\,F)$ is
the big tangent bundle of the manifold $M^{\mathcal{F}}$, which is
the set $M$ endowed with the differentiable structure of the sum
of the leaves of $\mathcal{F}$. Hence, $T^{big}\mathcal{F}$ has
the corresponding Courant structure (\ref{gC}), (\ref{crosetCou}).
A cross section of $T^{big}\mathcal{F}$ may be represented as
$Y\oplus\bar\alpha$ $(Y\in\chi(M^{\mathcal{F}}),\alpha\in
\Omega^1(M^{\mathcal{F}}))$, where the bar denotes the equivalence class
of $\alpha$ modulo $ann\,F$ (this bar-notation is always used
hereafter); generally, these cross sections are differentiable on
the sum of leaves. If we consider $Y_{l}\oplus\bar\alpha_{l}$
($l=1,2$) such that $Y_{l}\in\chi(M)$ and
$\alpha_{l}\in\Omega^1(M)$ are differentiable with respect to the
initial differentiable structure of $M$ we get the inner product
and Courant bracket
\begin{equation}\label{gF} g_F(Y_1\oplus\bar\alpha_1,
Y_2\oplus\bar\alpha_2)=\frac{1}{2}(\alpha_1(Y_2)+\alpha_2(Y_1)),\end{equation}
\begin{equation}\label{crosetF}
[Y_1\oplus\bar\alpha_1,Y_2\oplus\bar\alpha_2]
=([Y_1,Y_2]\oplus\overline{(L_{Y_1}\alpha_2-
L_{Y_2}\alpha_1+\frac{1}{2}d(\alpha_1(Y_2)-\alpha_2(Y_1))}),\end{equation}
where the results remain unchanged if $\alpha_{l}\mapsto
\alpha_{l}+\gamma_{l}$ with $\gamma_{l}\in ann\,F$. Formulas
(\ref{gF}), (\ref{crosetF}) show that $T^{big}\mathcal{F}\rightarrow
M$, where $M$ has its initial differentiable structure, is a Courant
algebroid with the anchor given by projection on the first term.
Alternatively, we can prove the same result by starting with
(\ref{gF}), (\ref{crosetF}) as definition formulas and by checking
the axioms of a Courant algebroid by computation.

We will transfer the Courant structure of $T^{big}\mathcal{F}$ by
the isomorphism $$ \Phi=id\oplus\hspace{1pt}^t\hspace{-1pt}\iota:
F\oplus (T^*M/ann\,F)\rightarrow F\oplus ann\,Q,$$ i.e.,
$$\Phi(Y\oplus\bar\alpha)=Y\oplus\alpha_{0,1},\;\;(Y\in F,\alpha=
\alpha_{1,0}+\alpha_{0,1}\in T^*M).$$ This makes $F\oplus ann\,Q$
into a Courant algebroid, which we shall denote by
$\mathcal{Q}=T^{big}_Q\mathcal{F}$, with the anchor equal to the
projection on $F$, the metric given by (\ref{gF}) and the bracket $$
[Y_1\oplus\bar\alpha_1, Y_2\oplus\bar\alpha_2]_{\mathcal{Q}}
=[Y_1,Y_2]\oplus pr_{ann\,Q}(L_{Y_1}\alpha_2-
L_{Y_2}\alpha_1+\frac{1}{2}d(\alpha_1(Y_2)-\alpha_2(Y_1))$$
$\alpha_1,\alpha_2\in ann\,Q$. Using the formula $L_Y=i(Y)d+di(Y)$
and the well known decomposition
$d=d'_{1,0}+d''_{0,1}+\partial_{2,-1}$ \cite{V0}, the expression of
the previous bracket becomes $$[Y_1\oplus\bar\alpha_1,
Y_2\oplus\bar\alpha_2]_{\mathcal{Q}}
=([Y_1,Y_2]\oplus(i(Y_1)d''\alpha_2-i(Y_2)d''\alpha_1$$
$$+\frac{1}{2}d''(i(Y_1)\alpha_2-i(Y_2)\alpha_1))\;\;(\alpha_1,\alpha_2\in
T^*_{0,1}M).$$
\section{The extension theorem}
Let $(M,\mathcal{F})$ be a foliated manifold. If the definition of a
Courant algebroid is modified by asking the anchor to be a morphism
$E\rightarrow\nu\mathcal{F}$, by asking $E,g,\sharp_E$ to be
foliated, by asking only for a bracket
$[\,,\,]_E:\Gamma_{fol}E\times\Gamma_{fol}E\rightarrow\Gamma_{fol}E$
and by asking the axioms to hold for foliated cross sections and
functions, then, we get the notion of a {\it transversal-Courant
algebroid} $(E,g_E,\sharp_E,[\,,\,]_E)$ over $(M,\mathcal{F})$
\cite{V1}. (The index $fol$ denotes foliated objects, i.e., objects
that either project to or are a lift of a corresponding object of
the space of leaves.)

On the other hand, a subbundle $B$ of a  Courant algebroid $A$
over $(M,\mathcal{F})$ is a {\it foliation} of $A$ if: i) $B$ is
$g_A$-isotropic and $\Gamma B$ is closed by $A$-brackets, ii)
$\sharp_A(B)=T\mathcal{F}$, iii) if $C=B^{\perp_{g_A}}$, then the
$A$-Courant structure induces the structure of a
transversal-Courant algebroid on the vector bundle $C/B$; then,
the pair $(A,B)$ is called a {\it foliated Courant algebroid} (see
\cite{V1} for details).

In this section we prove the announced result:\\
{\bf Theorem.} {\it Let $E$ be a transversal-Courant algebroid over
the foliated manifold $(M,\mathcal{F})$ and let $Q$ be a
complementary bundle of $F$ in $TM$. Then $E$ has a natural
extension to a foliated Courant algebroid $A$ with a foliation $B$
isomorphic to $F$.}
\begin{proof} The proof of this theorem requires a lot of technical
calculations. We will only sketch the path to be followed, leaving
the actual calculations to the interested reader. We shall denote
the natural extension that we wish to construct, and its operations,
by the index $0$. Take $A_0=T^{big}_Q\mathcal{F}\oplus
E=\mathcal{Q}\oplus E$ with the metric $g_0=g_F\oplus g_E$ and the
anchor $\sharp_0=pr_F\oplus\rho$, where $\rho=\varphi\circ\sharp_E$
with $\varphi$ defined by (\ref{descTM}), therefore,
$\psi\circ\rho=\sharp_E$. Notice that this implies
\begin{equation}\label{partial0}\partial_0\lambda=(0,\lambda|_F)+\frac{1}{2}
\sharp_{g_E}(\lambda\circ\rho)=(0,\lambda|_F)+\partial_E(\hspace{1pt}
^t\hspace{-1pt}\varphi\lambda)\;\;(\lambda\in T^*M)\end{equation}
and, in particular, $$\partial_0f =\partial_{\mathcal{Q}}(d''f)
\oplus\partial_E(d'f)=(0,d''f) \oplus\partial_E(d'f)\;\;(f\in
C^\infty(M)).$$

Then, inspired by the case $T^{big}M=
\mathcal{Q}\oplus\nu\mathcal{F}$ where the formulas below hold,
we define the bracket of
generating cross sections $Y\oplus\alpha\in\Gamma \mathcal{Q}$,
$e\in\Gamma_{fol}E$ by
\begin{equation}\label{cr0} \begin{array}{l}
[Y_1\oplus\alpha_1,Y_2,\oplus\alpha_2]_0=
[Y_1\oplus\alpha_1,Y_2,\oplus\alpha_2]_{\mathcal{Q}}\vspace{2mm}\\
\hspace*{1cm}\oplus\frac{1}{2}\sharp_{g_E}((L_{Y_1}\alpha_2-
L_{Y_2}\alpha_1+\frac{1}{2}d(\alpha_1(Y_2)-\alpha_2(Y_1))\circ\rho)\vspace{2mm}\\

=([Y_1,Y_2]\oplus0)+\partial_0(L_{Y_1}\alpha_2-
L_{Y_2}\alpha_1+\frac{1}{2}d(\alpha_1(Y_2)-\alpha_2(Y_1)),\vspace{2mm}\\

[e,Y\oplus\alpha]_0= ([\rho e,Y]\oplus (L_{\rho e}\alpha)|_F)
\oplus\frac{1}{2}\sharp_{g_E}((L_{\rho e}\alpha)\circ\rho)\vspace{2mm}\\

\hspace*{1cm}=([\rho e,Y]\oplus0)+\partial_0L_{\rho e}\alpha),\vspace{2mm}\\

[e_1,e_2]_0=(([\rho e_1,\rho e_2]-\rho [e_1,e_2]_E)\oplus0) \oplus
[e_1,e_2]_E.\end{array}\end{equation} The first term of the right
hand side of the second formula belongs to $\Gamma \mathcal{Q}$
since $e\in\Gamma_{fol}E$ implies $[\rho e,Y]\in\Gamma F$. The
first term of the right hand side of the third formula belongs to
$\Gamma
\mathcal{Q}$ since we have
$$\psi([\rho e_1,\rho e_2]-\rho [e_1,e_2]_E)=\psi([\rho e_1,\rho e_2])-
\sharp_E [e_1,e_2]_E$$ $$=\psi([\rho e_1,\rho e_2])-
[\sharp_Ee_1,\sharp_Ee_2]_{\nu\mathcal{F}}=0.$$

Furthermore, we extend the bracket (\ref{cr0}) to arbitrary cross
sections in agreement with the axiom 4) of Courant algebroids,
i.e., for any functions $f,f_1,f_2\in C^\infty(M)$, we define
\begin{equation}\label{f} \begin{array}{l} [Y\oplus\alpha,fe]_0=
f[Y\oplus\alpha,e]_0\oplus(Yf)e,\vspace{2mm}\\

[f_1e_1,f_2e_2]_0=f_1f_2[e_1,e_2]_0+ f_1(\rho e_1(f_2))e_2
-f_2(\rho e_2(f_1))e_1\vspace{2mm}\\

\hspace*{2cm}-g_E(e_1,e_2)(f_1\partial_0f_2-f_2\partial_0f_1)\end{array}
\end{equation}
($Y\in\Gamma F,\alpha\in ann\,Q, e,e_1,e_2\in\Gamma_{fol}E$). It
follows easily that formulas (\ref{cr0}) and  (\ref{f}) give the
same result if $f\in C^\infty_{fol}(M,\mathcal{F})$.

We have to check that the bracket defined by (\ref{cr0}),
(\ref{f}) satisfies the axioms of a Courant algebroid and it is
enough to do that for every possible combination of arguments of
the form $Y\oplus\alpha\in
\mathcal{Q}$ and $fe$, $e\in\Gamma_{fol}E$, $f\in
C^\infty(M)$.

To check axiom 1), apply the anchor $\sharp_0=pr_F+\rho$ to each
of the five formulas (\ref{cr0}), (\ref{f}) and use the
transversal-Courant algebroid axioms satisfied by $E$. To check
axiom 2), use formula (\ref{partial0}). The required results
follow straightforwardly. It is also easy to check axiom 4) from
(\ref{f}) and from axiom 4) for $\mathcal{Q}$ and $E$.

Furthermore, technical (lengthy) calculations show that if we have
a bracket such that axioms 1), 2), 4) hold, then, if 5) holds for
a triple of arguments, 5) also holds if the same arguments are
multiplied by arbitrary functions. Therefore, in our case it
suffices to check axiom 5) for the following six triples: (i)
$(e,e_1,e_2)$, (ii) $(Y\oplus\alpha,e_1,e_2)$, (iii)
$(e,Y\oplus\alpha,e')$, (iv) $(Y\oplus\alpha,Y'\oplus\alpha',e)$,
(v) $(e,Y_1\oplus\alpha_1,Y_2\oplus\alpha_2)$, (vi)
$(Y\oplus\alpha,Y_1\oplus\alpha_1,Y_2\oplus\alpha_2)$, where all
$Y\oplus\alpha\in\Gamma\mathcal{Q}$ and all $e\in\Gamma_{fol}E$.
In cases (i), (vi) the result follows from axiom 5) satisfied by
$E,\mathcal{Q}$, respectively. In the other cases computations
involving evaluations of Lie derivatives will do the job.

Finally, we have to check axiom 3). If we consider any vector
bundle $E$ with an anchor and a bracket that satisfy axioms 1),
2), 4), 5), then, by applying axiom 5) to the triple $(e,
e_1=\partial f,e_2)$ $(f\in C^\infty(M))$ we get
\begin{equation}\label{crosetpt5} [e,\partial
f]_E=\frac{1}{2}\partial(\sharp_Ee(f)),\end{equation} whence (using
local coordinates, for instance) the following general formula
follows
\begin{equation}\label{gen-partial} [e,\partial_E\alpha]_E=
\partial_E(L_{\sharp_Ee}\alpha-\frac{1}{2}d(\alpha(\sharp_Ee))).
\end{equation}

Furthermore, assuming again that axioms 1), 2), 4), 5) hold and
using (\ref{partialcug}) and (\ref{crosetpt5}) a lengthy but
technical calculation shows that, if axiom 3) holds for a triple
$(e_1,e_2,e_3)$, it also holds for $(e_1,e_2,fe_3)$ $(f\in
C^\infty(M))$ provided that
\begin{equation}\label{Erond} \begin{array}{c} \mathcal{E}:=
g([e_1,e_2],e_3)+\frac{1}{2}\sharp e_2(g(e_1,e_3))-
\frac{1}{2}\sharp
e_1(g(e_2,e_3))\vspace{2mm}\\=\frac{1}{3}\sum_{Cycl}g([e_1,e_2],e_3)\end{array}
\end{equation} ($:=$ denotes a definition). But, if the last two
terms in $\mathcal{E}$ are expressed by axiom 5) for $E$ followed
by (\ref{crosetpt5}), and after we repeat the same procedure one
more time, we get
$$ \mathcal{E}=\frac{1}{4}\sum_{Cycl}g([e_1,e_2],e_3)
+\frac{1}{4}\mathcal{E},$$ whence we see that (\ref{Erond}) holds
for any triple $(e_1,e_2,e_3)$.

Hence, it suffices to check axiom 3) for the following cases: (i)
$(e_1,e_2,e_3)$, (ii) $(e_1,e_2,Y\oplus\alpha)$, (iii)
$(e,Y_1\oplus\alpha_1,Y_2\oplus\alpha_2)$, (iv)
$(Y_1\oplus\alpha_1,Y_2\oplus\alpha_2,Y_3\oplus\alpha_3)$, where all
$Y\oplus\alpha\in\Gamma\mathcal{Q}$ and all $e\in\Gamma_{fol}E$. In
case (i), using the second and third formula (\ref{cr0}), we get
$$[[e_1,e_2]_0,e_3]_0=([[\rho e_1,\rho e_2],\rho e_3]-
\rho[[e_1,e_2]_E,e]_E,0)\oplus[[e_1,e_2]_E,e]_E$$ and the required
result follows in view of the Jacobi identity for vector fields and
of axiom 3) for $E$ (in this case, the right hand side of axiom 3)
for $A_0$ reduces to the one for $E$).

To check the result in the other cases simpler, we decompose
$$Y\oplus\alpha=(Y\oplus0)+(0\oplus\alpha)$$ and check axiom 3) for
each case induced by this decomposition.

For a triple $(e_1,e_2,Y\oplus0)$ the right hand side of axiom 3)
is zero and the left-hand side is $$([[\rho e_1,\rho
e_2],Y]+[[\rho e_2,Y],\rho e_1]+[[Y,\rho e_1],\rho
e_2])\oplus0=0$$ by the Jacobi identity for vector fields.

For a triple $(e_1,e_2,0\oplus\alpha)$, after cancelations, the
right hand side of axiom 3) becomes $(1/2)\alpha([\rho e_1,\rho
e_2])$. The same result is obtained for the left hand side if we use
the second form of the first two brackets defined by (\ref{cr0}) and
formula (\ref{gen-partial}).

For a triple $(e,Y_1\oplus0,Y_2\oplus0)$ the two sides of axiom 3)
vanish (the left hand side reduces to the Jacobi identity for the
vector fields $(\rho e,Y_1,Y_2)$), hence the axiom holds.

For a triple $(e,0\oplus\alpha_1,0\oplus\alpha_2)$, using the second
form of the second bracket (\ref{cr0}) and formula
(\ref{gen-partial}), axiom 3) reduces to $0=0$, i.e., the axiom
holds.

For a triple $(e,Y_1\oplus0,0\oplus\alpha)$, if we notice that
$0\oplus\alpha=\partial_0\alpha$ (see (\ref{partial0})) and use
(\ref{gen-partial}), we see that the two sides of the equality
required by axiom 3) are equal to $(1/2)\partial_0(\alpha([\rho
e,Y])-(1/2)\rho e(\alpha(Y)))$, hence the axiom holds.

The case $(Y_1\oplus0,Y_2\oplus0,Y_3\oplus0)$ is trivial. In the
case $(Y_1\oplus0,Y_2\oplus0,0\oplus\alpha=\partial_0\alpha)$
similar computations give the value
$(1/4)\partial_0(\alpha([Y_1,Y_2])-d\alpha(Y_1,Y_2))$ for the two
sides of the corresponding expression of axiom 3). Finally, in the
cases $(Y\oplus0,\partial_0\alpha_1,\partial_0\alpha_2)$ and
$(\partial_0\alpha_1,\partial_0\alpha_2,\partial_0\alpha_3)$ the two
sides of the required equality are $0$ since the image of
$\partial_0$ is isotropic and the restriction of the bracket to this
image is zero (use axiom 2) and formula (\ref{gen-partial})).
\end{proof}
\hspace*{7.5cm}{\small \begin{tabular}{l} Department of
Mathematics\\ University of Haifa, Israel\\ E-mail:
vaisman@math.haifa.ac.il \end{tabular}}
\end{document}